\documentclass{amsart}
\usepackage{amsfonts}
\usepackage{amscd}
\usepackage{amssymb}
\usepackage{mathrsfs}
\theoremstyle{plain}
\newtheorem{thm}{Theorem}[section]
\newtheorem{cor}[thm]{Corollary}
\newtheorem{lem}[thm]{Lemma}
\newtheorem{pro}[thm]{Proposition}
\theoremstyle{definition}

\newtheorem{defn}[thm]{Definition}
\usepackage{amsmath}
\usepackage{amssymb}
\usepackage{amsthm}
\usepackage{latexsym}
\usepackage[T1]{fontenc}
\usepackage[all]{xy}
\DeclareMathOperator{\Per}{Per} \DeclareMathOperator{\supp}{supp}

\begin{document}
\author{Christian Svensson}
\address{Mathematical Institute, Leiden University,
P.O. Box 9512, 2300 RA Leiden, The Netherlands, and Centre for
Mathematical Sciences, Lund University, Box 118, SE-221 00 Lund,
Sweden} \email{chriss@math.leidenuniv.nl}
\author{Sergei Silvestrov}
\address{Centre for Mathematical Sciences,
Lund University, Box 118, SE-221 00 Lund, Sweden}
\email{Sergei.Silvestrov@math.lth.se}
\author{Marcel de Jeu}
\address{Mathematical Institute,
Leiden University, P.O. Box 9512, 2300 RA Leiden, The Netherlands}
\email{mdejeu@math.leidenuniv.nl}

\title[Connections between dynamical systems
and crossed products ...]
{Connections between dynamical systems and \\
crossed products of Banach algebras by
$\mathbb{Z}$.}



\noindent \keywords{ Crossed product; Banach
algebra; ideals, dynamical system; maximal
abelian subalgebra \\ \noindent {\it
\textup{2000} Mathematics Subject
Classification.} Primary 47L65 Secondary 16S35,
37B05, 54H20} \noindent
\date{January 30, 2007}

\begin{abstract}
Starting with a complex commutative semi-simple regular Banach
algebra $A$ and an automorphism $\sigma$ of $A$, we form the crossed
product of $A$ with the integers, where the latter act on $A$ via
iterations of $\sigma$. The automorphism induces a topological
dynamical system on the character space $\Delta(A)$ of $A$ in a
natural way. We prove an equivalence between the property that every
non-zero ideal in the crossed product has non-zero intersection with
the subalgebra $A$, maximal commutativity of $A$ in the crossed
product, and density of the non-periodic points of the induced
system on the character space. We also prove that every non-trivial
ideal in the crossed product always intersects the commutant of $A$
non-trivially. Furthermore, under the assumption that $A$ is unital
and such that $\Delta(A)$ consists of infinitely many points, we
show equivalence between simplicity of the crossed product and
minimality of the induced system, and between primeness of the
crossed product and topological transitivity of the system.
\end{abstract}
\maketitle
\section{Introduction}
A lot of work has been done in the direction of
connections between certain topological dynamical
systems and crossed product $C^*$-algebras. In
\cite{tom1} and \cite{tom2}, for example, one
starts with a homeomorphism $\sigma$ of a compact
Hausdorff space $X$ and constructs the crossed
product $C^*$-algebra $C(X) \rtimes_{\alpha}
\mathbb{Z}$, where $C(X)$ is the algebra of
continuous complex valued functions on $X$ and
$\alpha$ is the automorphism of $C(X)$ naturally
induced by $\sigma$. One of many results obtained
is equivalence between simplicity of the algebra
and minimality of the system, provided that $X$
consists of infinitely many points, see
\cite{davidson}, \cite{power}, \cite{tom1},
\cite{tom2} or, for a more general result,
\cite{williams}. In \cite{svesildej}, a purely
algebraic variant of the crossed product is
considered, and with more general classes of
algebras than merely continuous functions on
compact Hausdorff spaces serving as ``coefficient
algebras'' in the crossed products. For example,
it was proved there that, for such crossed
products, the analogue of the equivalence between
density of non-periodic points of a dynamical
system and maximal commutativity of the
``coefficient algebra'' in the associated crossed
product \mbox{$C^*$-algebra} is true for
significantly larger classes of coefficient
algebras and associated dynamical systems. In
this paper, we go beyond these results and
investigate the ideal structure of some of the
crossed products considered in \cite{svesildej}.
More specifically, we consider crossed products
of complex commutative semi-simple regular Banach
algebras $A$ with the integers under an
automorphism $\sigma : A \rightarrow A$.

In Section~\ref{defbas} we give the most general definition of the
kind of crossed product that we will use throughout this paper. We
also mention the elementary result that the commutant of the
coefficient algebra is automatically a commutative subalgebra of the
crossed product. The more specific setup that we will be working in
is introduced in Section~\ref{setbas}. There we also introduce some
notation and mention two basic results concerning a canonical
isomorphism between certain crossed products, and an explicit
description of the commutant of the coefficient algebra in one of
them.

According to \cite[Theorem 5.4]{tom2}, the following three
properties are equivalent:
\begin{itemize}
\item The non-periodic points of $(X, \sigma)$ are dense in $X$;
\item Every non-zero closed ideal $I$ of the crossed product
 $C^*$-algebra $C(X) \rtimes_{\alpha} \mathbb{Z}$ is such that $I \cap C(X) \neq
 \{0\}$;
\item $C(X)$ is a maximal abelian $C^*$-subalgebra of $C(X)
 \rtimes_{\alpha} \mathbb{Z}$.
\end{itemize}
In Section~\ref{treeq}, an analogue of this result is proved for our
setup. A reader familiar with the theory of crossed product
$C^*$-algebras will easily recognize that if one chooses $A = C(X)$
for $X$ a compact Hausdorff space in Corollary~\ref{triquiv} below,
then the crossed product is canonically isomorphic to a norm-dense
subalgebra of the crossed product $C^*$-algebra coming from the
considered induced dynamical system. We also combine this with a
theorem from \cite{svesildej} to conclude a stronger result for the
Banach algebra $L_1 (G)$, where $G$ is a locally compact abelian
group with connected dual group.

In Section~\ref{simplmin}, we prove the equivalence between
algebraic simplicity of the crossed product and minimality of the
induced dynamical system in the case when $A$ is unital with its
character space consisting of infinitely many points. This is
analogous to \cite[Theorem 5.3]{tom2}, \cite[Theorem VIII
3.9]{davidson}, the main result in \cite{power} and, as a special
case of a more general result, \cite[Corollary 8.22]{williams} for
the crossed product $C^*$-algebra.

In Section~\ref{comint}, the fact that the commutant of $A$ always
has non-zero intersection with any non-zero ideal of the crossed
product is shown. This should be compared with the fact that $A$
itself may well have zero intersection with such ideals, as
Corollary~\ref{triquiv} shows. The analogue of this result in the
context of crossed product $C^*$-algebras appears to be open.

Finally, in Section~\ref{primtra} we show equivalence between
primeness of the crossed product and topological transitivity of the
induced system, in the case when $A$ is unital and has an infinite
character space. The analogue of this in the context of crossed
product $C^*$-algebras is~\cite[Theorem 5.5]{tom2}.
\section{Definition and a basic result}\label{defbas}
Let $A$ be an associative commutative complex algebra and let $\Psi
: A \rightarrow A$ be an algebra automorphism. Consider the set
\[A \rtimes_{\Psi} \mathbb{Z} = \{f: \mathbb{Z} \rightarrow A \,|\,
 f(n) = 0 \,\,\textup{except for a finite number of}\, \,n\}.\]
We endow it with the structure of an associative complex algebra by
defining scalar multiplication and addition as the usual pointwise
operations. Multiplication is defined by \emph{twisted convolution},
$*$, as follows;
\[(f*g) (n) = \sum_{k \in \mathbb{Z}} f(k) \cdot \Psi^k (g(n-k)),\]
where $\Psi^k$ denotes the $k$-fold composition of $\Psi$ with
itself. It is trivially verified that $A \rtimes_{\sigma}
\mathbb{Z}$ \emph{is} an associative $\mathbb{C}$-algebra under
these operations. We call it the \emph{crossed product of $A$ and
$\mathbb{Z}$ under $\Psi$}.

A useful way of working with $A \rtimes_{\Psi}
\mathbb{Z}$ is to write elements $f, g \in A
\rtimes_{\Psi} \mathbb{Z}$ in the form $f =
\sum_{n \in \mathbb{Z}} f_n \delta^n$, $g =
\sum_{m \in \mathbb{Z}} g_m \delta^m$, where $f_n
= f(n)$, $g_m = g(m)$, addition and scalar
multiplication are canonically defined, and
multiplication is determined by $(f_n
\delta^n)*(g_m \delta^m) = f_n \cdot \Psi^n (g_m)
\delta^{n+m}$, where $n,m \in \mathbb{Z}$ and
$f_n, g_m \in A$ are arbitrary.

Clearly one may canonically view $A$ as an abelian subalgebra of $A
\rtimes_{\Psi} \mathbb{Z}$, namely as $\{f_0 \delta^0 \, | \, f_0
\in A\}$. The following elementary result is proved in
\cite[Proposition 2.1]{svesildej}.

\begin{pro}\label{comcom}
The commutant $A'$ of $A$ in $A \rtimes_{\Psi}
\mathbb{Z}$ is abelian, and thus it is the unique
maximal abelian subalgebra containing $A$.
\end{pro}
\section{Setup and two basic results}\label{setbas}
In what follows, we shall focus on cases when $A$
is a commutative complex Banach algebra, and
freely make use of the basic theory for such $A$,
see e.g. \cite{larsen}. As conventions tend to
differ slightly in the literature, however, we
mention that we call a commutative Banach algebra
$A$ \emph{semi-simple} if the Gelfand transform
on $A$ is injective, and that we call it
\emph{regular} if, for every subset $F $ of the
character space $\Delta(A)$ that is closed in the
Gelfand topology and for every $\phi_0 \in
\Delta(A) \setminus F$, there exists an $a \in A$
such that $\widehat{a} (\phi) = 0$ for all $\phi
\in F$ and $\widehat{a} (\phi_0) \neq 0$. All
topological considerations of the character space
$\Delta(A)$ will be done with respect to its
Gelfand topology (the weakest topology making all
elements in the image of the Gelfand transform of
$A$ continuous on $\Delta(A)$).

Now let $A$ be a complex commutative semi-simple regular Banach
algebra, and let $\sigma  : A \rightarrow A$ be an algebra
automorphism. As in \cite{svesildej}, $\sigma$ induces a map
$\widetilde{\sigma}: \Delta(A) \rightarrow \Delta(A)$ (where
$\Delta(A)$ denotes the character space of $A$) defined by
$\widetilde{\sigma}(\mu) = \mu \circ \sigma^{-1}, \mu \in
\Delta(A)$, which is automatically a homeomorphism when $\Delta(A)$
is endowed with the Gelfand topology. Hence we obtain a topological
dynamical system $(\Delta(A), \widetilde{\sigma})$. In turn,
$\widetilde{\sigma}$ induces an automorphism $\widehat{\sigma} :
\widehat{A} \rightarrow \widehat{A}$ (where $\widehat{A}$ denotes
the algebra of Gelfand transforms of all elements of $A$) defined by
$\widehat{\sigma} (\widehat{a}) = \widehat{a} \circ
\widetilde{\sigma}^{-1} = \widehat{\sigma(a)}$. Therefore we can
form the crossed product $\widehat{A} \rtimes_{\widehat{\sigma}}
\mathbb{Z}$. We also mention that when speaking of ideals, we will
always mean two-sided ideals.

In what follows, we shall make frequent use of the following fact.
Its proof consists of a trivial direct verification.
\begin{thm}\label{isom}
Let $A$ be a commutative semi-simple Banach
algebra and $\sigma$ be an automorphism, inducing
an automorphism $\widehat{\sigma} : \widehat{A}
\rightarrow \widehat{A}$ as above. Then the map
$\Phi : A \rtimes_{\sigma} \mathbb{Z} \rightarrow
\widehat{A} \rtimes_{\widehat{\sigma}}
\mathbb{Z}$ defined by $\sum_{n \in \mathbb{Z}}
a_n \delta^n \mapsto \sum_{n \in \mathbb{Z}}
\widehat{a_n} \delta^n$ is an isomorphism of
algebras mapping $A$ onto $\widehat{A}$.
\end{thm}
Before stating the next result,
we make the following basic
definitions.
\begin{defn}
For any nonzero $n \in \mathbb{Z}$ we set
\begin{align*}\label{per}
\Per^n (\Delta(A)) &= \{\mu \in \Delta(A)\,|\, \mu =
 \widetilde{\sigma}^n (\mu)\}. \\
\end{align*}
Furthermore, we denote the non-periodic points by
\begin{align*}
\Per^{\infty}(\Delta(A)) &= \bigcap_{n \in \mathbb{Z} \setminus \{0\}}
 (\Delta(A) \backslash \Per^n (\Delta(A)).\\
\end{align*}
Finally, for $f \in \widehat{A}$, put
\begin{align*}\supp(f) &= \overline{\{\mu \in \Delta(A) \, | \, f(\mu)
 \neq 0\}}.
\end{align*}
\end{defn}

\begin{thm}\label{commdesc} We have the following
explicit description of ${\widehat{A}}^{\,'}$ in
$\widehat{A} \rtimes_{\widehat{\sigma}}
\mathbb{Z}$:
\[{\widehat{A}}^{\,'} =\{\sum_{n \in \mathbb{Z}} f_n \delta^n \,|\, f_n
 \in \widehat{A}, \,\textup{and for all} \, n \in \mathbb{Z}:
 \supp(f_n) \subseteq \Per^n (\Delta(A))\}.\]
\end{thm}
\noindent {\em Proof.} This follows from
\cite[Corollary 3.4]{svesildej}, as $\widehat{A}$
trivially separates the points of $\Delta(A)$ and
$\Per^n (\Delta(A))$ is a closed set. \qed

\section{Three equivalent properties}\label{treeq}
In this section we shall conclude that, for certain $A$, two
different algebraic properties of $A \rtimes_{\sigma} \mathbb{Z}$
are equivalent to density of the non-periodic points of the
naturally associated dynamical system on the character space
$\Delta(A)$, and hence obtain equivalence of three different
properties. The analogue of this result in the context of crossed
product $C^*$-algebras is \cite[Theorem 5.4]{tom2}. We shall also
combine this with a theorem from \cite{svesildej} to conclude a
stronger result for the Banach algebra $L_1 (G)$, where $G$ is a
locally compact abelian group with connected dual group.

In \cite[Theorem 4.8]{svesildej}, the following result is proved.
\begin{thm}\label{eqcomperden}
Let $A$ be a complex commutative regular semi-simple Banach algebra,
$\sigma: A \rightarrow A$ an automorphism and $\widetilde{\sigma}$
the homeomorphism of $\Delta(A)$ in the Gelfand topology induced by
$\sigma$ as described above. Then the non-periodic points are dense
in $\Delta(A)$ if and only if $\widehat{A}$ is a maximal abelian
subalgebra of $\widehat{A} \rtimes_{\widehat{\sigma}} \mathbb{Z}$.
In particular, $A$ is maximal abelian in $A \rtimes_{\sigma}
\mathbb{Z}$ if and only if the non-periodic points are dense in
$\Delta(A)$.
\end{thm}
We shall soon prove another algebraic property of the crossed
product equivalent to density of the non-periodic points of the
induced system on the character space. First, however, we need two
easy topological lemmas.

\begin{lem}\label{sepophd}
Let $x \in \Delta(A)$ be such that the points $\widetilde{\sigma}^i
(x)$ are distinct for all $i$ such that $-m\leq i \leq n$, where $n$
and $m$ are positive integers. Then there exist an open set $U_x$
containing $x$ such that the sets $\widetilde{\sigma}^{i} (U_x)$ are
pairwise disjoint for all $i$ such that $-m \leq i \leq n$.
\end{lem}
\noindent{\em Proof.} It is easily checked that
any finite set of points in a Hausdorff space can
be separated by pairwise disjoint open sets.
Separate the points $\widetilde{\sigma}^i (x)$
with disjoint open sets $V_i$. Then it is readily
verified that the set \[U_x :=
\widetilde{\sigma}^{m} (V_{-m}) \cap
\widetilde{\sigma}^{m-1} (V_{-m+1}) \cap \ldots
\cap V_0 \cap \widetilde{\sigma}^{-1} (V_1) \cap
\ldots \cap \widetilde{\sigma}^{-n} (V_n)\] is an
open neighbourhood of $x$ with the required
property. \qed
\begin{lem}\label{bairegrej}
The non-periodic points of $(\Delta(A),
\widetilde{\sigma})$ are dense if and only if the
set $\Per^{n} (\Delta(A))$ has empty interior for
all positive integers $n$.
\end{lem}
\noindent {\em Proof.} Clearly, if there is a
positive integer $n_0$ such that $\Per^{n_0}
(\Delta(A))$ has non-empty interior, the
non-periodic points are not dense. For the
converse, we recall that $\Delta(A)$ is a Baire
space since it is locally compact and Hausdorff,
and note that we may write
\[\Delta(A) \backslash \Per^{\infty} (\Delta(A)) = \bigcup_{n > 0}
 \Per^n (\Delta(A)).\]
If the set of non-periodic points is not dense,
its complement has non-empty interior, and as the
sets $\Per^n (\Delta(A))$ are clearly all closed,
there must exist an integer $n_0 > 0$ such that
$\Per^{n_0} (\Delta(A))$ has non-empty interior
since $\Delta(A)$ is a Baire space. \qed

We are now ready to prove the promised result.
\begin{thm}\label{aperidov}
Let $A$ be a complex commutative semi-simple regular Banach algebra,
$\sigma: A \rightarrow A$ an automorphism and $\widetilde{\sigma}$
the homeomorphism of $\Delta(A)$ in the Gelfand topology induced by
$\sigma$ as described above. Then the non-periodic points are dense
in $\Delta(A)$ if and only if every non-zero ideal $I \subseteq A
\rtimes_{\sigma} \mathbb{Z}$ is such that $I \cap A \neq \{0\}$.
\end{thm}
\noindent {\em Proof.} We first assume that
$\overline{\Per^\infty (\Delta(A))} = \Delta(A)$,
and work initially in $\widehat{A}
\rtimes_{\widehat{\sigma}} \mathbb{Z}$. Assume
that $I \subseteq \widehat{A}
\rtimes_{\widehat{\sigma}} \mathbb{Z}$ is a
non-zero ideal, and that $f = \sum_{n \in
\mathbb{Z}} f_n \delta^n \in I$. By definition,
only finitely many $f_n$ are non-zero. Denote the
set of integers $n$ for which $f_n \not \equiv 0$
by $S=\{n_1, \ldots, n_r\}$. Pick a non-periodic
point $x \in \Delta(A)$ such that $f_{n_1} (x)
\neq 0$ (by density of $\Per^\infty (\Delta(A))$
such $x$ exists). Using the fact that $x$ is not
periodic we may, by Lemma~\ref{sepophd}, choose
an open neighbourhood $U_x$ of $x$ such that
$\widetilde{\sigma}^{-n_i} (U_x) \cap
\widetilde{\sigma}^{-n_j}(U_x) = \emptyset$ for
$n_i \neq n_j,\, n_i, n_j \in S$. Now by
regularity of $A$ we can find a function $g \in
\widehat{A}$ that is non-zero in
$\widetilde{\sigma}^{-n_1} (x)$, and vanishes
outside $\widetilde{\sigma}^{-n_1}(U_x)$.
Consider $f*g = \sum_{n \in \mathbb{Z}} f_n \cdot
(g \circ \widetilde{\sigma}^{-n}) \delta^n$. This
is an element in $I$ and clearly the coefficient
of $\delta^{n_1}$ is the only one that does not
vanish on the open set $U_x$. Again by regularity
of $A$, there is an $h \in \widehat{A}$ that is
non-zero in $x$ and vanishes outside $U_x$.
Clearly $h*f*g = [h \cdot (g \circ
\widetilde{\sigma}^{-n_1}) f_{n_1}] \delta^{n_1}$
is a non-zero monomial belonging to $I$. Now any
ideal that contains a non-zero monomial
automatically contains a non-zero element of
$\widehat{A}$. Namely, if $a_i \delta^i \in I$
then $[a_i \delta^i] * [(a_i \circ
\widetilde{\sigma}^{i}) \delta^{-i}] = a_i ^2 \in
\widehat{A}$. By the canonical isomorphism in
Theorem~\ref{isom}, the result holds for $A
\rtimes_{\sigma} \mathbb{Z}$ as well.

For the converse, assume that $\overline{\Per^{\infty} (\Delta(A))}
\neq \Delta(A)$. Again we work in $\widehat{A}
\rtimes_{\widehat{\sigma}} \mathbb{Z}$. It follows from
Lemma~\ref{bairegrej} that since $\overline{\Per^{\infty}
(\Delta(A))} \neq \Delta(A)$, there exists an integer $n >0$ such
that $\Per^n (\Delta(A))$ has non-empty interior. As $A$ is assumed
to be regular, there exists $f \in \widehat{A}$ such that $\supp(f)
\subseteq \Per^n (\Delta(A))$. Consider now the ideal $I=(f + f
\delta^n)$. It  consists of finite sums of elements of the form $a_i
\delta^i (f + f \delta^n) a_j \delta^j$. Using that $f$ vanishes
outside $\Per^n (\Delta(A))$, we may rewrite this as follows
\begin{align*}
&a_i \delta^i (f + f \delta^n) a_j \delta^j = [a_i \cdot (a_j \circ
\widetilde{\sigma}^{-i})\delta^{i}]* [f \delta^j + f \delta^{n+j}]\\
&= [a_i \cdot (a_j \circ \widetilde{\sigma}^{-i}) \cdot (f\circ
\widetilde{\sigma}^{-i})] \delta^{i+j} + [a_i \cdot (a_j \circ
\widetilde{\sigma}^{-i}) \cdot (f\circ \widetilde{\sigma}^{-i})]
\delta^{i+j+n}.
\end{align*}
This means that any element in $I$ may be written
in the form $\sum_{i} (b_i \delta^i +
b_i\delta^{n+i})$. As $i$ runs only through a
finite subset of $\mathbb{Z}$, this is not a
non-zero monomial. In particular, it is not a
non-zero element in $\widehat{A}$. Hence $I$
intersects $\widehat{A}$ trivially. By the
canonical isomorphism in Theorem~\ref{isom}, the
result carries over to $A \rtimes_{\sigma}
\mathbb{Z}$. \qed

Combining
Theorem~\ref{eqcomperden} and
Theorem~\ref{aperidov}, we now have the following
result.
\begin{cor}\label{triquiv}
Let $A$ be a complex commutative semi-simple regular Banach algebra,
$\sigma: A \rightarrow A$ an automorphism and $\widetilde{\sigma}$
the homeomorphism of $\Delta(A)$ in the Gelfand topology induced by
$\sigma$ as described above. Then the following three properties are
equivalent:
\begin{itemize}
\item The non-periodic points $\Per^\infty (\Delta(A))$ of
 $(\Delta(A),\widetilde{\sigma})$ are dense in $\Delta(A)$;
\item Every non-zero ideal $I \subseteq A \rtimes_{\sigma} \mathbb{Z}$
 is such that $I \cap A \neq \{0\}$;
\item $A$ is a maximal abelian subalgebra of $A \rtimes_{\sigma}
 \mathbb{Z}$.
\end{itemize}
\end{cor}
We shall make use of Corollary~\ref{triquiv} to conclude a result
for a more specific class of Banach algebras. We start by recalling
a number of standard results from the theory of Fourier analysis on
groups, and refer to \cite{larsen} and \cite{rudin} for details. Let
$G$ be a locally compact abelian group. Recall that $L_1 (G)$
consists of equivalence classes of complex valued Borel measurable
functions of $G$ that are integrable with respect to a Haar measure
on $G$, and that $L_1 (G)$ equipped with convolution product is a
commutative regular semi-simple Banach algebra. A group homomorphism
$\gamma : G \rightarrow \mathbb{T}$ from a locally compact abelian
group to the unit circle is called a \emph{character} of $G$. The
set of all \emph{continuous} characters of $G$ forms a group
$\Gamma$, the \emph{dual group} of $G$, if the group operation is
defined by
\[(\gamma_1 + \gamma_2) (x) = \gamma_1 (x) \gamma_2 (x) \,\, (x \in G;
 \gamma_1, \gamma_2 \in \Gamma).\]
If $\gamma \in \Gamma$ and if
\[\widehat{f}(\gamma) = \int_{G} f(x) \gamma (-x) dx \quad (f \in L_1
 (G)),\]
then the map $f \mapsto \widehat{f} (\gamma)$ is a non-zero complex
homomorphism of $L_1 (G)$. Conversely, every non-zero complex
homomorphism of $L_1 (G)$ is obtained in this way, and distinct
characters induce distinct homomorpisms. Thus we may identify
$\Gamma$ with $\Delta(L_1 (G))$. The function $\widehat{f} : \Gamma
\rightarrow \mathbb{C}$ defined as above is called the \emph{Fourier
transform} of $f \in L_1 (G)$, and is precisely the Gelfand
transform of $f$. We denote the set of all such $\widehat{f}$ by
$A(\Gamma)$. Furthermore, $\Gamma$ is a locally compact abelian
group in the Gelfand topology.

In \cite[Theorem 4.16]{svesildej}, the following result is proved.
\begin{thm}\label{conndual}
Let $G$ be a locally compact abelian group with connected dual group
and let $\sigma : L_1 (G) \rightarrow L_1 (G)$ be an automorphism.
Then $L_1 (G)$ is maximal abelian in $L_1 (G) \rtimes_{\sigma}
\mathbb{Z}$ if and only if $\sigma$ is not of finite order.
\end{thm}
Combining Corollary~\ref{triquiv} and Theorem~\ref{conndual} the
following result is immediate.
\begin{cor}\label{conndualtri}
Let $G$ be a locally compact abelian group with connected dual group
and let $\sigma : L_1 (G) \rightarrow L_1 (G)$ be an automorphism.
Then the following three statements are equivalent.
\begin{itemize}
\item $\sigma$ is not of finite order;
\item Every non-zero ideal $I \subseteq L_1 (G) \rtimes_{\sigma}
 \mathbb{Z}$ is such that $I \cap L_1 (G) \neq \{0\}$;
\item $L_1 (G)$ is a maximal abelian subalgebra of $L_1 (G)
 \rtimes_{\sigma} \mathbb{Z}$.
\end{itemize}
\end{cor}
\section{Minimality versus simplicity}\label{simplmin}
Recall that a topological dynamical system is said to be
\emph{minimal} if all of its orbits are dense, and that an algebra
is called \emph{simple} if it lacks non-trivial proper ideals.
\begin{thm}\label{simpmin}
Let $A$ be a complex commutative semi-simple regular unital Banach
algebra such that $\Delta(A)$ consists of infinitely many points,
and let $\sigma : A \rightarrow A$ be an algebra automorphism of
$A$. Then $A \rtimes_{\sigma} \mathbb{Z}$ is simple if and only if
the naturally induced system $(\Delta(A), \widetilde{\sigma})$ is
minimal.
\end{thm}
\noindent {\em Proof.} Suppose first that the
system is minimal, and assume that $I$ is a
proper ideal of $A \rtimes_{\sigma} \mathbb{Z}$.
Note that $I \cap A$ is a proper $\sigma$- and
$\sigma^{-1}$-invariant ideal of $A$. By basic
theory of Banach algebras, $I \cap A$ is
contained in a maximal ideal of $A$ (note that $I
\cap A \neq A$ as $A$ is unital and $I$ was
assumed to be proper), which is the kernel of an
element $\mu \in \Delta(A)$. Now $\widehat{I \cap
A}$ is a $\widehat{\sigma}$- and
$\widehat{\sigma}^{-1}$-invariant proper
non-trivial ideal of $\widehat{A}$, all of whose
elements vanish in $\mu$. Invariance of this
ideal implies that all of its elements even
annihilate the whole orbit of $\mu$ under
$\widetilde{\sigma}$. But by minimality, every
such orbit is dense and hence $\widehat{I \cap A}
= \{0\}$. By semi-simplicity of $A$, this means
$I \cap A = \{0\}$, so $I=\{0\}$ by
Corollary~\ref{triquiv}. For the converse, assume
that there is an element $\mu \in \Delta(A)$
whose orbit $O(\mu)$ is not dense. By regularity
of $A$ there is a nonzero $g \in \widehat{A}$
that vanishes on $\overline{O(\mu)}$. Then
clearly the ideal generated by $g$ in
$\widehat{A} \rtimes_{\widehat{\sigma}}
\mathbb{Z}$ consists of finite sums of elements
of the form $(f_n \delta^n) * g * (h_m \delta^m)
= [f_n \cdot(g \circ \widetilde{\sigma}^{-n})
\cdot(h_m \circ
\widetilde{\sigma}^{-n})]\delta^{n+m}$, and hence
the coefficient of every power of $\delta$ in
this ideal must vanish in $\mu$, whence the ideal
is proper. Hence by Theorem~\ref{isom}, $A
\rtimes_{\sigma} \mathbb{Z}$ is not simple. \qed

\section{Every non-zero ideal has non-zero intersection with
 $A'$}\label{comint}
We shall now show that any non-zero ideal in $A
\rtimes_{\sigma} \mathbb{Z}$ has non-zero
intersection with $A'$. This should be compared
with Corollary~\ref{triquiv}, which says that a
non-zero ideal may well intersect $A$ solely in
$0$. We have no information on the validity of
the analogue of this result in the context of
crossed product $C^*$-algebras.

\begin{thm}\label{commint}
Let $A$ be a complex commutative semi-simple
regular Banach algebra, and $\sigma : A
\rightarrow A$ an automorphism. Then every
non-zero ideal $I$ in $A \rtimes_{\sigma}
\mathbb{Z}$ has non-zero intersection with the
commutant $A'$ of $A$ in $A \rtimes_{\Psi}
\mathbb{Z}$, that is\, $I \cap A' \neq \{0\}$.
\end{thm}
\noindent {\em Proof.} As usual, we work in
$\widehat{A} \rtimes_{\widehat{\sigma}}
\mathbb{Z}$. When $\overline{\Per^{\infty}
(\Delta(A))} = \Delta(A)$, the result follows
immediately from Corollary~\ref{triquiv}. We will
use induction on the number of non-zero terms in
an element $f= \sum_{n \in \mathbb{Z}} f_n
\delta^n$ to show that it generates an ideal that
intersects $\widehat{A}^{\,'}$ non-trivially. The
starting point for the induction, namely when $f
= f_n \delta^n$ with non-zero $f_n$, is clear
since any such element generates an ideal that
even intersects $\widehat{A}$ non-trivially, as
was shown in the proof of Theorem~\ref{aperidov}.
Now assume inductively that the conclusion of the
theorem is true for the ideals generated by any
element of $\widehat{A}
\rtimes_{\widehat{\sigma}} \mathbb{Z}$ with $r$
non-zero terms for some positive integer $r$, and
consider an element $f = f_{n_1} \delta^{n_1} +
\ldots + f_{n_{r+1}} \delta^{n_{r+1}}$. By
multiplying from the right with a suitable
element we obtain an element in the ideal
generated by $f$ of the form $g =
\sum_{i=0}^{m_r} g_{i}\delta^i$ such that $g_0
\not \equiv 0$. If some of the other $g_i$ are
zero we are done by induction hypothesis, so we
may assume this is not the case. We may also
assume that $g$ is not in the commutant of
$\widehat{A}$ since otherwise we are of course
also done. This means, by Theorem~\ref{commdesc},
that there is such $j$ that $0 < j \leq m_r$ and
$\supp(g_j) \not \subseteq \Per^j (\Delta(A))$.
Pick an $x \in \supp(g_j)$ such that $x \neq
\widetilde{\sigma}^{-j} (x)$ and $g_j (x) \neq
0$. As $\Delta(A)$ is Hausdorff we can choose an
open neighbourhood $U_x$ of $x$ such that $U_x
\cap \widetilde{\sigma}^{-j} (U_x) = \emptyset$.
Regularity of $A$ implies existence of an $h \in
\widehat{A}$ such that $h \circ
\widetilde{\sigma}^{-j} (x) =1$ and $h$ vanishes
identically outside of $\widetilde{\sigma}^{-j}
(U_x)$. Now $g*h = \sum_{i=0}^{m_r}  g_i \cdot (h
\circ \widetilde{\sigma}^{-i}) \delta^i$. Using
regularity of $A$ again we pick a function $a \in
\widehat{A}$ such that $a (x) =1$ and $a$
vanishes outside $U_x$. We have $a*g*h =
\sum_{i=0}^{m_r}  a \cdot g_i \cdot (h \circ
\widetilde{\sigma}^{-i}) \delta^i$, which is in
the ideal generated by $f$. Now $a \cdot g_0
\cdot h$ is identically zero since $a \cdot h =
0$. On the other hand, $a \cdot g_j \cdot (h\circ
\widetilde{\sigma}^{-j})$ is non-zero in the
point $x$. Hence $a * g
* h$ is a non-zero element in the ideal generated by $f$ whose
number of non-zero coefficient functions is less
than or equal to $r$. By the induction
hypothesis, such an element generates an ideal
that intersects the commutant of $\widehat{A}$
non-trivially. By Theorem~\ref{isom} it follows
that every non-zero ideal in
$A\rtimes_{\sigma}\mathbb{Z}$ intersects $A'$
non-trivially. \qed
\section{Primeness versus topological transitivity}\label{primtra}
We shall show that for certain $A$, $A \rtimes_{\sigma} \mathbb{Z}$
is prime if and only if the induced system $(\Delta(A),
\widetilde{\sigma})$ is topologically transitive. The analogue of
this result in the context of crossed product $C^*$-algebras is
in~\cite[Theorem 5.5]{tom2}.
\begin{defn}
The system $(\Delta(A), \widetilde{\sigma})$ is called
\emph{topologically transitive} if for any pair of non-empty open
sets $U, V$ of $\Delta(A)$, there exists an integer $n$ such that
$\widetilde{\sigma}^n (U) \cap V \neq \emptyset$.
\end{defn}
\begin{defn}
The algebra $A \rtimes_{\sigma} \mathbb{Z}$ is called \emph{prime}
if the intersection between any two non-zero ideals $I, J$ is
non-zero, that is $I \cap J \neq \{0\}$.
\end{defn}
For convenience, we also make the following definition.
\begin{defn}
The map $E : \widehat{A} \rtimes_{\widehat{\sigma}} \mathbb{Z}
\rightarrow \widehat{A}$ is defined by $E(\sum_{n \in \mathbb{Z}}
f_n \delta^n) = f_0$.
\end{defn}
To prove the main theorem of this section, we need the two following
topological lemmas.
\begin{lem}\label{toptraset}
If $(\Delta(A), \widetilde{\sigma})$ is not topologically
transitive, then there exist two disjoint invariant non-empty open
sets $O_1$ and $O_2$ such that $\overline{O_1} \cup \overline{O_2} =
\Delta(A)$.
\end{lem}
\noindent {\em Proof.} As the system is not
topologically transitive, there exist non-empty
open sets $U, V \subseteq \Delta(A)$ such that
for any integer $n$ we have $\widetilde{\sigma}^n
(U) \cap V = \emptyset$. Now clearly the set $O_1
= \bigcup_{n \in \mathbb{Z}} \widetilde{\sigma}^n
(U)$ is an invariant non-empty open set. Then
$\overline{O_1}$ is an invariant closed set. It
follows that $O_2 = \Delta(A) \backslash
\overline{O_1}$ is an invariant open set
containing $V$. Thus we even have that
$\overline{O_1} \cup O_2 = \Delta(A)$, and the
result follows. \qed
\begin{lem}\label{toptraper}
If $(\Delta(A), \widetilde{\sigma})$ is topologically transitive and
there is an $n_0 > 0$ such that $\Delta(A) = \Per^{n_0}
(\Delta(A))$, then $\Delta(A)$ consists of a single orbit and is
thus finite.
\end{lem}
\noindent {\em Proof.} Assume two points $x, y
\in \Delta(A)$ are not in the same orbit. As
$\Delta(A)$ is Hausdorff we may separate the
points $x, \sigma (x), \ldots, \sigma^{n_0 -1}
(x), y$ by pairwise disjoint open sets $V_0, V_1,
\ldots, V_{n_0 -1}, V_y$. Now consider the set
\[U_x := V_0 \cap \widetilde{\sigma}^{-1} (V_1) \cap
 \widetilde{\sigma}^{-2} (V_2) \cap \ldots \cap \widetilde{\sigma}^{-n_0 +1} (V_{n_0
 -1}).\]
Clearly the sets $A_x  = \cup_{i=0}^{n_0 -1}
\widetilde{\sigma}^i (U_x)$ and $A_y =
\cup_{i=0}^{n_0 -1} \widetilde{\sigma}^i (V_y)$
are disjoint invariant non-empty open sets, which
leads us to a contradiction. Hence $\Delta(A)$
consists of one single orbit under
$\widetilde{\sigma}$. \qed

We are now ready for a proof of the following
result.
\begin{thm}\label{pritop}
Let $A$ be a complex commutative semi-simple regular unital Banach
algebra such that $\Delta(A)$ consists of infinitely many points,
and let $\sigma$ be an automorphism of $A$. Then $A \rtimes_{\sigma}
\mathbb{Z}$ is prime if and only if the associated system
$(\Delta(A), \widetilde{\sigma})$ on the character space is
topologically transitive.
\end{thm}
\noindent {\em Proof.} Suppose $(\Delta(A),
\widetilde{\sigma})$ is not topologically
transitive. Then there exists, by
Lemma~\ref{toptraset}, two disjoint invariant
non-empty open sets $O_1$ and $O_2$ such that
$\overline{O_1} \cup \overline{O_2} = \Delta(A)$.
Let $I_1$ and $I_2$ be the ideals generated in
$\widehat{A} \rtimes_{\widehat{\sigma}}
\mathbb{Z}$ by $k (\overline{O_1})$ (the set of
all functions in $\widehat{A}$ that vanish on
$\overline{O_1}$) and $k (\overline{O_2})$
respectively. We have that
\[E(I_1 \cap I_2) \subseteq E(I_1) \cap E(I_2) = k(\overline{O_1}) \cap
 k(\overline{O_2})= k(\overline{O_1} \cup \overline{O_2})= k(\Delta(A))
 = \{0\}.\]
It is not difficult to see that if $I \subseteq
\widehat{A} \rtimes_{\widehat{\sigma}}
\mathbb{Z}$ is an ideal and $E(I) = \{0\}$, then
$I=\{0\}$. Namely, suppose $F= \sum_{n} f_n
\delta^n \in I$ and $f_i \neq 0$ for some integer
$i$. Since $A$ is unital, so is $\widehat{A}$ and
thus $\delta^{-1}\in \widehat{A}
\rtimes_{\widehat{\sigma}} \mathbb{Z}$. So $F*
\delta^{-i} \in I$ and hence $E(F* \delta^{-i}) =
f_i = 0$ which is a contradiction, so $I =
\{0\}$. Hence $I_1 \cap I_2 = \{0\}$ and
$\widehat{A} \rtimes_{\widehat{\sigma}}
\mathbb{Z}$ is not prime. By Theorem~\ref{isom},
neither is $A \rtimes_{\sigma} \mathbb{Z}$. Next
suppose that $(\Delta(A), \widetilde{\sigma})$ is
topologically transitive. Assume that
$\Per^\infty (\Delta(A))$ is not dense. Then by
Lemma~\ref{bairegrej} there is an integer $n_0 >
0$ such that $\Per^{n_0} (\Delta(A))$ has
non-empty interior. As $\Per^{n_0} (\Delta(A))$
is invariant and closed, topological transitivity
implies that $\Delta(A) = \Per^{n_0}
(\Delta(A))$. This, however, is impossible since
by Lemma~\ref{toptraper} it would force
$\Delta(A)$ to consist of a single orbit and
hence be finite. Thus $\Per^\infty (\Delta(A))$
is dense after all. Now let $I$ and $J$ be two
non-zero proper ideals in $A \rtimes_{\sigma}
\mathbb{Z}$. Unitality of $A$ assures us that $I
\cap A$ and $J \cap A$ are proper invariant
ideals of $A$ and density of $\Per^\infty
(\Delta(A))$ assures us that they are non-zero,
by Theorem~\ref{triquiv}. Consider $A_I = \{\mu
\in \Delta(A) \,|\, \mu (a) = 0 \, \textup{for
all} \,\,  a \in I \cap A\}$ and $A_J = \{\nu \in
\Delta(A) \,|\, \nu (b) = 0 \, \textup{for all}
\,\, b \in J \cap A\}$. Now by Banach algebra
theory a proper ideal in a commutative unital
Banach algebra $A$ is contained in a maximal
ideal, and a maximal ideal of $A$ is always
precisely the set of zeroes of some $\xi \in
\Delta(A)$. This implies that both $A_I$ and
$A_J$ are non-empty, and semi-simplicity of $A$
assures us that they are proper subsets of
$\Delta(A)$. Clearly they are also closed and
invariant under $\widetilde{\sigma}$ and
$\widetilde{\sigma}^{-1}$. Hence $\Delta(A)
\setminus A_I$ and $\Delta(A) \setminus A_J$ are
invariant non-empty open sets. By topological
transitivity we must have that these two sets
intersect, hence that $A_I \cup A_J \neq
\Delta(A)$. This means that there exists $\eta
\in \Delta(A)$ and $a \in I \cap A$, $b \in J
\cap A$ such that $\eta(a) \neq 0$, $\eta(b) \neq
0$ and hence that $\eta (ab) \neq 0$. Hence $0
\neq ab \in I \cap J$, and we conclude that $A
\rtimes_{\sigma} \mathbb{Z}$ is prime. \qed

\subsection*{Acknowledgments} This work was supported by a visitor's
grant of the \textit{Netherlands Organisation for
Scientific Research (NWO)}, \textit{The Swedish
Foundation for International Cooperation in
Research and Higher Education (STINT)},
\textit{Crafoord Foundation, The Royal
Physiographic Society in Lund, The Royal Swedish
Academy of Sciences and The Swedish Research
Council}.

We are also grateful to Jun Tomiyama and Johan {\"O}inert for useful
discussions.

\end{document}